\documentclass[12pt]{article}
\usepackage{amssymb,latexsym,amsmath}
\usepackage{amsfonts}
\usepackage[cp1251]{inputenc}
\usepackage[epsf]{}
\usepackage[english,russian]{babel}
\usepackage[dvips]{graphicx}
\newtheorem{theorem}{Theorem}

\newtheorem{proposition}{Proposition}

\begin{document}
\title{Supporting manifolds for high-dimensional Morse-Smale diffeomorphisms with few saddles}
\author{Medvedev V.\and Zhuzhoma E.}
\date{}
\maketitle

\renewcommand{\figurename}{Figure}
\renewcommand{\abstractname}{Abstract}
\renewcommand{\refname}{Bibliography}

\begin{abstract}
We describe the topological structure of closed manifolds of dimension no less than four which admit Morse-Smale diffeomorphisms such that its non-wandering set contains any number of sink periodic points, and any number of source periodic points, and few saddle periodic points.
\end{abstract}

\section*{Introduction}

In 1960, Smale \cite{Smale60a} introduced the class of dynamical systems (flows and diffeomorphisms) later called Morse-Smale systems. It was proved that the Morse-Smale systems are structurally stable ones with zero topological entropy \cite{Palis69,PalisSmale70,Robinson99}. In this sense, the Morse-Smale systems are simplest structurally stable dynamical systems. Originally, the Morse-Smale systems was introduced as dynamical systems whose non-wandering sets consist of a finite number of hyperbolic orbits with transversal inter\-sec\-tions of invariant manifolds \cite{Smale60a,Smale67}.
There are deep connections between dynamics and topological structures of supporting ma\-ni\-folds \cite{GrinesGurevichZhPochinka2019,MedvedevZhuzhoma2008-MIAN,Smale61a}. On a closed manifold, any Morse-Smale system has at least one attracting orbit and at least one repelling orbit \cite{Smale60a}. Thus, the simplest Morse-Smale diffeomorphism has the non-wandering set consisting of two points: a sink and source. In this case, the supporting $n$-manifold is an $n$-sphere denoted by $\mathbb{S}^n$, and such orientation preserving Morse-Smale diffeomorphism $\mathbb{S}^n\to\mathbb{S}^n$ is a so called north-south (or, source-sink) diffeomorphism \cite{GrinesMedvedevPochinkaZh2010,Reeb1952}.
In \cite{MedvedevZhuzhoma2013-top-appl}, the authors described completely the topological structure of closed supporting manifolds for the Morse-Smale diffeomorphisms whose wandering sets consist of three points,
see Proposition \ref{prop:from-Medvedev-Zh-Top-Appl-2013}. It is natural to study the topological structure of closed manifolds admitting the Morse-Smale diffeomorphisms whose non-wandering sets consist of $\geq 4$ points. Obviously, if we replace an original diffeomorphism by any its iteration, the topological structure of supporting manifold does not change. Therefore, without loss of generality, one can consider the Morse-Smale diffeomorphisms whose non-wandering sets consist of fixed points. Later on, $M^n$ is a closed smooth connected $n$-manifold.

Denote by $MS(M^n;a,b,c)$ the set of Morse-Smale diffeomorphisms $f: M^n\to M^n$ such that the non-wandering set of $f$ consists of $a$ sinks, $b$ sources, and $c$ saddles. In addition, we assume that the restriction of $f$ on any invariant manifold of every saddle is a preserving orientation mapping.
Recall that due to Smale \cite{Smale60a}, $a\geq 1$ and $b\geq 1$, and if $c=0$ then $M^n=\mathbb{S}^n$. Below, we suggest $c\geq 1$. There are numerous papers devoted to the periodic dates of Morse-Smale systems
for the dimensions $n=1,2,3$, see the articles \cite{Batterson79-81,BattersonHandelNarasimhan81,BlanchardFranks80,Narasimhan79} and the surveys \cite{GrinesGurevichZhPochinka2019,MedvedevZhuzhoma2008-MIAN}.
Therefore, we suppose $n\geq 4$. Mainly, we pay attention for $c=1,2$, and $c\geq 3$ codimension one saddles. Ideologically, this paper is a continuation of the article \cite{MedvedevZhuzhoma2013-top-appl} where the authors considered the set $MS(M^n;1,1,1)$.

First, we begin with the description of supporting manifolds and non-wandering sets for the Morse-Smale diffeomorp\-hisms $MS(M^n;a,b,1)$ with a unique saddle.
Recall that the Morse index of fixed point $\sigma$ equals the dimension of the unstable manifold of $\sigma$. Later on, $\mathbb{B}^n$ is an open $n$-ball, and $\mathbb{S}^k$ is a $k$-sphere. The union $A\sqcup B$ means a disjoint union i.e., $A\cap B=\emptyset$.
\begin{theorem}\label{thm:one-saddle-structure}
Let $MS(M^n;a,b,1)$ be the set of Morse-Smale diffeomorphisms of closed smooth $n$-manifold $M^n$, $n\geq 4$, with a unique saddle denoted by $\sigma$. Then either

1) $a=b=1$; moreover,

1a) $n\in\{4,8,16\}$;

1b) $M^n=\mathbb{B}^n\sqcup\mathbb{S}^{\frac{n}{2}}$;

1c) the saddle $\sigma$ has $\frac{n}{2}$-dimensional separatrices;

1d) the homotopy groups $\pi_1(M^n)=\cdots=\pi_{\frac{n}{2}-1}(M^n)=0$, and hence, $M^n$ is a simply connected orientable manifold;

or

2) $a+b=3$, i.e. either $a=1$, $b=2$ or $a=2$, $b=1$; moreover,

2a) $M^n=\mathbb{S}^n$,

2b) the Morse index of $\sigma$ equals $(n-1)$ provided $a=1$, $b=2$, and the Morse index of $\sigma$ equals 1 provided $a=2$, $b=1$.
\end{theorem}

Recall that a codimension one saddle has two one-dimensional separatrices and one $(n-1)$-dimensional separatrix.
Denote by $N\otimes S^1$ the total space of locally trivial fiber bundle $N\otimes S^1\longrightarrow S^1$ with the base circle $S^1$ and fiber a manifold $N$. Let $\mathbb{D}^k$ be a closed $k$-disk. The following statements concern to the class $MS(M^n;a,b,2)$ with two saddles.
\begin{theorem}\label{thm:two-saddles-at-least-one-codimension-one}
Let $MS(M^n;a,b,2)$ be the set of Morse-Smale diffeomorphisms of closed smooth $n$-manifold $M^n$, $n\geq 4$, with two saddles denoted by $\sigma_1$, $\sigma_2$. Then

I) If $\sigma_1$, $\sigma_2$ are saddles of codimension one each then either

1.1) $a=b=1$, and $M^n$ is the union of two copies $\mathbb{D}^{n-1}\otimes S^1$; in addition, the saddles $\sigma_1$, $\sigma_2$ have different Morse indexes,

or

1.2) $a+b=4$ (i.e., either $a=b=2$, or $a=1$, $b=3$, or $a=3$, $b=1$), and $M^n=\mathbb{S}^n$; in addition,

1.2a) $(a,b)=(2,2)$ iff the saddles have different Morse indexes;

1.2b) $(a,b)=(1,3)$ iff the saddles both have the Morse index equals $n-1$;

1.2c) $(a,b)=(3,1)$ iff the saddles both have the Morse index equals 1.

II) If $\sigma_1$ is a codimension one saddle while $\sigma_2$ is not a codimension one and codimension two saddle, then

2.1) $a+b=3$, i.e. either $a=1$, $b=2$ or $a=2$, $b=1$; moreover

2.2) $n\in\{4,8,16\}$, and $M^n=\mathbb{B}^n\sqcup\mathbb{S}^{\frac{n}{2}}$;

2.3) the homotopy groups $\pi_1(M^n)=\cdots=\pi_{\frac{n}{2}-1}(M^n)=0$, and hence, $M^n$ is a simply connected and orientable manifold.
\end{theorem}

The next statement concerns to $c\geq 3$ saddles of codimension one.
\begin{theorem}\label{thm-all-saddles-codimension-one-same-morse-index}
Let $MS(M^n;a,b,c)$ be the set of Morse-Smale diffeomorphisms of closed smooth $n$-manifold $M^n$, $n\geq 4$, with $c\geq 3$ saddles. If all saddles has the same Morse index $\in\{1;n-1\}$, then
$M^n=\mathbb{S}^n$ and $(a,b)\in \{(1,c+1),(c+1,1)\}$. In addition,

1) $(a,b)=(1,c+1)$ iff the Morse index of all saddles equals $n-1$;

2) $(a,b)=(c+1,1)$ iff the Morse index of all saddles equals 1.
\end{theorem}

For the Morse-Smale diffeomorphisms without codimension one saddles at all, we prove the following result.
\begin{theorem}\label{thm:there-are-no-cod1-at-all}
Let $MS(M^n;a,b,c)$ be the set of Morse-Smale diffeomorphisms of closed smooth $n$-manifold $M^n$, $n\geq 4$, with no codimension one saddles. Then $a=b=1$, and $M^n$ is a simply connected manifold.
\end{theorem}

%***********************************************************************************************

The structure of the paper is the following. In Section \ref{s:d-and-prev}, we formulate the main definitions and give some previous results. In Section \ref{s:proofs}, one proves main results.
In Section \ref{s:exm-and-conclusions}, one constructs examples showing that the sets of Morse-Smale diffeomorphisms under consideration are not empty.

%************************************************************************************************

\textit{Acknowledgments}. The authors are partially supported by Laboratory of Dynamical Systems and Applications of National Research University Higher School of Economics, of the Ministry of science and higher education of the RF, grant ag. № 075-15-2019-1931, except Theorem 2 supported by RNF, grant 17-11-01041.

\section{Basic definitions and previous results}\label{s:d-and-prev}

Basic definitions of Dynamical Systems see in \cite{ABZ}, \cite{Robinson99}, \cite{Smale67}.
Let $f: M^n\to M^n$ be a diffeomorphism of $n$-manifold $M^n$, and $p$ a periodic point of period $k\in\mathbb{N}$. The stable manifold $W^s(p)$ is defined to be the set of points $x\in M^n$ such that $\varrho(f^{kj}(x);p)\to 0$ as $j\to\infty$ where $\varrho$ is a metric on $M^n$. The unstable manifold $W^u(p)$ is the stable manifold of $p$ for the diffeomorphism $f^{-1}$. Stable and unstable manifolds are called invariant manifolds. It is well known that if $p$ is hyperbolic, then every invariant manifold is an immersed submanifold homeomorphic to Euclidean space. Moreover, $W^s(p)$ and $W^u(p)$ are intersected transversally at $p$, and $\dim W^s(p)+\dim W^u(p)=n$.

A diffeomor\-p\-hism $f$ is \textit{Morse-Smale} if it is structurally stable and the non-wandering set $NW(f)$ of $f$ consists of a finitely many periodic orbits. In particular, each periodic orbit is hyperbolic and, stable and unstable manifolds of periodic orbits intersect transversally.

Let $f: M^n\to M^n$ be a Morse-Smale diffeomorphism of $n$-manifold $M^n$.
A periodic orbit $p$ is called a \textit{sink periodic point} (resp. \textit{source periodic point}) if $\dim W^s(p)=n$ and $\dim W^u(p)=0$ (resp. $\dim W^s(p)=0$ and $\dim W^u(p)=n$) where $W^u(p)$ and $W^s(p)$ are unstable and stable manifolds of $p$ respectively. A sink or source periodic orbit is called a \textit{node periodic orbit}.
A periodic point $\sigma$ is called a \textit{saddle periodic point} if
$1\leq\dim W^u(\sigma)\leq n-1$, $1\leq\dim W^s(\sigma)\leq n-1$ where $W^u(\sigma)$ and $W^s(\sigma)$ are unstable and stable manifolds of $\sigma$ respectively.
A component of $W^u(\sigma)\setminus\sigma$ denoted by $Sep^u(\sigma)$ is called
an \textit{unstable separatrix} of $\sigma$. If $\dim W^u(\sigma)\geq 2$, then $Sep^u(\sigma)$ is a unique component.

A saddle periodic point $\sigma$ is a \textit{codimension one} saddle periodic point if either $\dim W^u(\sigma)=1$, $\dim W^s(\sigma)=n-1$ or $\dim W^u(\sigma)=n-1$, $\dim W^s(\sigma)=1$.
If $\dim W^u(\sigma)=1$, then there are two one-dimensional separatrices denoted by $Sep^u_1(\sigma)$ and $Sep^u_2(\sigma)$. The similar notation holds for a stable separatrix.

A sink (resp., source, saddle) periodic point $p$ is called a \textit{sink} (resp., \textit{source, saddle}) if $p$ is a fixed point.

Let $p$, $q$ are saddle periodic points such that $W^u(p)\cap W^s(q)\neq\emptyset$. The intersection $W^u(p)\cap W^s(q)$ is called \textit{heteroclinic}. Due to the transversality $W^u(p)\pitchfork W^s(q)$, a heteroclinic intersection is either points or $m$-submanifolds with $m\geq 1$. If $\dim\left(W^u(p)\cap W^s(q)\right)\geq 1$, the intersection $W^u(p)\cap W^s(q)$ is called a \textit{heteroclinic manifold}.
For the reference, we formulate the following result proved in \cite{GinesMedvedevZh2003} (see also \cite{GrinesMedvedevPochinkaZh2010}).
\begin{proposition}\label{prop:closure-is-sphere-for-dim-more-one}
Let $f: M^n\to M^n$ be a Morse-Smale diffeomorphism, and $W^{\tau}(\sigma)$ an invariant manifold of dimension $2\leq d\leq n-1$ of a saddle $\sigma$ where $\tau\in\{s,u\}$. Suppose that $W^{\tau}(\sigma)$ has no heteroclinic intersections with other separatrices. Then $Sep^{\tau}(\sigma)=W^{\tau}(\sigma)\setminus\{\sigma\}$ belongs to unstable (if $\tau=s$) or stable (if $\tau=u$) manifold of some node (source or sink, respectively) periodic point, say $N$, and the topological closure of $Sep^{\tau}(\sigma)$ is a topologically embedded $d$-sphere that equals $W^{\tau}(\sigma)\cup\{N\}=S^d_{\tau}$. Moreover, $S^d_{\tau}$ is a locally flat embedded sphere provided
$d\neq n-2$.
\end{proposition}

For $1\le m\le n$, we presume Euclidean space $\mathbb{R}^m$ to be included naturally in $\mathbb{R}^n$ as the subset whose final $(n-m)$ coordinates each equals $0$. Let $e: M^m\to N^n$ be an embedding of a boundaryless $m$-manifold $M^m$ in the interior of $n$-manifold $N^n$. One says that $e(M^m)$ is \textit{locally flat at} $e(x)$, $x\in M^m$, if there exists a neighborhood $U(e(x))=U$ and a homeomorphism $h: U\to\mathbb{R}^n$ such that
$ h(U\cap e(M^m))=\mathbb{R}^m\subset\mathbb{R}^n. $ Otherwise, $e(M^m)$ is \textit{wild at} $e(x)$. One says that $e(M^m)$ is a \textit{locally flat embedded} submanifold if $e(M^m)$ is locally flat at every point $e(x)$,
$x\in M^m$ \cite{DavermanVenema-book-2009}. The similar notation holds for a compact $M^m$ with a boundary, in particular $M^m=[0;1]$.

Note that a separatrix $Sep^{\tau}(\sigma)$ is a smooth manifold. Hence, $Sep^{\tau}(\sigma)$ is locally flat at every point \cite{DavermanVenema-book-2009}. However a-priori, the topological closure
$clos Sep^{\tau}(\sigma)=W^{\tau}(\sigma)\cup\{N\}$ could be wild at nodes of this topological closure.

Now, we prove the statement we'll need bellow. The end of proof will be denoted by $\diamondsuit$.
\begin{proposition}\label{prop:closure-is-sphere-for-dim-equals-one}
Let $f: M^n\to M^n$ be a Morse-Smale diffeomorphism, and $W^{\tau}(\sigma)$ an invariant manifold of dimension $d=1$ of a saddle $\sigma$ where $\tau\in\{s,u\}$. Suppose that $W^{\tau}(\sigma)$ has no heteroclinic intersections with other separatrices. Then every separatrix $Sep^{\tau}_i(\sigma)=W^{\tau}(\sigma)\setminus\{\sigma\}$, $i-1,2$, belongs to unstable (if $\tau=s$) or stable (if $\tau=u$) manifold of some node (source or sink, respectively) periodic point, say $N_i$, and the topological closure of $Sep^{\tau}_1(\sigma)\cup Sep^{\tau}_2(\sigma)$ is either a locally flat embedded circle or a locally flat embedded segment. In addition,

1) if $clos~\left(Sep_1(\sigma)\cup Sep_2(\sigma)\right)$ is a circle $S^1_0$ then $S^1_0=W^{\tau}(\sigma)\cup\{N\}$ where $n_0=N_1=N_2$, and $S_0$ has a neighborhood $T$ homeomorphic to $\mathbb{D}^{n-1}\otimes S^1$ such that $T$ is forward-invariant provided $n_0$ is a sink and back-invariant provided $n_0$ is a source. Moreover, $T$ contains only two fixed points: the saddle $\sigma$ and a node $n_0$.

2) if $clos~\left(Sep_1(\sigma)\cup Sep_2(\sigma)\right)$ is a segment $I$ then $I=Sep_1(\sigma)\cup Sep_2(\sigma)\cup\{\sigma\}\cup\{n_1\}\cup\{n_2\}$ where $n_1$, $n_2$ are nodes, and $I$ has a neighborhood $B$ homeomorphic to an $n$-ball such that $B$ is forward-invariant provided $n_1$, $n_2$ are sinks and back-invariant provided $n_1$, $n_2$ are sources. Moreover, $B$ contains only three fixed points: the saddle $\sigma$ and two nodes $n_1$, $n_2$.
\end{proposition}
\textsl{Proof of Proposition \ref{prop:closure-is-sphere-for-dim-equals-one}}. Because of $W^{\tau}(\sigma)$ has no heteroclinic intersections, the topological closure of $Sep^{\tau}_1(\sigma)\cup Sep^{\tau}_2(\sigma)$ is either a topologically embedded circle or a topologically embedded segment. Since $n\geq 4$, the codimension of $Sep^{\tau}_1(\sigma)\cup Sep^{\tau}_2(\sigma)$ is more than two. It follows from \cite{DavermanVenema-book-2009} that
$clos \left(Sep_1(\sigma)\cup Sep_2(\sigma)\right)$ is either a locally flat embedded circle or a locally flat embedded segment.

Suppose that $clos~\left(Sep_1(\sigma)\cup Sep_2(\sigma)\right)=S^1_0$ is a locally flat embedded circle. Then $S_0$ has a tubular neighborhood homeomorphic to $\mathbb{D}^{n-1}\otimes S^1$. Suppose for definiteness that $Sep_1(\sigma)$, $Sep_2(\sigma)$ are unstable separatrices. Then $S_0$ is the union $Sep_1(\sigma)\cup Sep_2(\sigma)\cup\{\sigma\}\cup\{n_0\}$ where $n_0$ is a sink. Taking $T$ smaller if necessary, one can assume that $T$ contains only two fixed points, $sigma$ and $n_0$. According to the Grobman-Hartman theorem, $T$ can choose to be forward-invariant i.e. $f(T)\subset T$. When $Sep_1(\sigma)$, $Sep_2(\sigma)$ are stable separatrices, the proof is similar.

Suppose that $clos~\left(Sep_1(\sigma)\cup Sep_2(\sigma)\right)=I$ is a locally flat embedded segment. Suppose for definiteness that $Sep_1(\sigma)$, $Sep_2(\sigma)$ are unstable separatrices. Then $I$ is the union
$Sep_1(\sigma)\cup Sep_2(\sigma)\cup\{\sigma\}\cup\{n_1\}\cup\{n_2\}$ where $n_1$, $n_2$ are sinks. Due to the Grobman-Hartman theorem, $I$ has a neighborhood $B$ homeomorphic to an $n$-ball, and $B$ can choose to be forward-invariant i.e. $f(B)\subset B$. Hence, $B$ contains only three fixed points $\sigma$, $n_1$, $n_2$. When $Sep_1(\sigma)$, $Sep_2(\sigma)$ are stable separatrices, the proof is similar.
$\diamondsuit$

We say that a Morse-Smale diffeomorphism $f: M^n\to M^n$ is \textit{without heteroclinic manifolds on codimension one separatrices} if given any codimension one saddle periodic point $p\in NW(f)$, the codimension one separatrix of $p$ does not contain heteroclinic manifolds. We'll need the following statement proved in \cite{GinesMedvedevZh2017}.
\begin{proposition}\label{prop:rough-decomposition}
Let $M^n$ be a closed $n$-manifold, $n\ge 3$, supporting a Morse-Smale diffeomorphism $f$ without heteroclinic manifolds on codimension one separatrices. Suppose that the non-wandering set $NW(f)$ consists of $\mu$ node periodic points, $\nu$ codimension one saddle periodic points, and arbitrary number of saddle periodic points that are not codimension one.
Then the number
 $$ g=\frac{1}{2}\left(\nu - \mu +2\right)\geq 0 $$
is integer. In addition, if $g=0$ then $M^n$ is either $\mathbb{S}^n$ or
\begin{equation}\label{eq:decomposition-1}
    M^n=N^n_1\sharp\cdots\sharp\, N^n_{l}
\end{equation}
for some $1\leq l\leq 1+\nu$ where every $N^n_i$ admits a polar Morse-Smale diffeomorphism without codimension one saddle periodic orbit.
\end{proposition}

For references, we quote one more result proved in \cite{MedvedevZhuzhoma2013-top-appl}.
\begin{proposition}\label{prop:from-Medvedev-Zh-Top-Appl-2013}
Suppose $f\in MS(M^n;1,1,1)$, $n\geq 2$, has a unique saddle $\sigma$. Then
\begin{itemize}
  \item $n\in\{2,4,8,16\}$;
  \item $M^n\mathbb{B}^n\sqcup \mathbb{S}^{\frac{n}{2}}$;
  \item $M^2$ is the projective plane;
  \item for $n\geq 4$, the homotopy groups $\pi_1(M^n)=\cdots=\pi_{\frac{n}{2}-1}(M^n)=0$, and hence, $M^n$ is simply connected and orientable;
  \item the stable and unstable manifolds $W^s(\sigma)$, $W^u(\sigma)$ respectively are both $\frac{n}{2}$-dimensional;
  \item  $n\in\{2,8,16\}$, $W^s(\sigma)\cup\{\alpha\}$ and $W^u(\sigma)\cup\{\omega\}$ are locally flat $\frac{n}{2}$-spheres where $\omega$ is a sink, $\alpha$ is a source;
  \item $M^n=W^s(\omega)\sqcup\left(W^s(\sigma)\cup\{\alpha\}\right)=W^u(\alpha)\sqcup\left(W^u(\sigma)\cup\{\omega\}\right)$.
\end{itemize}
\end{proposition}

\section{Proof of main results}\label{s:proofs}

A cutting of $M^n$ along a codimension one submanifold $N^{n-1}\subset M^n$ means that we delete a sufficiently small neighborhood $U$ of $N^{n-1}$ homeomorphic to $N^{n-1}\times (0;1)$, and take $clos(M^n\setminus U)$.

\medskip
\textsl{Proof of Theorem \ref{thm:one-saddle-structure}}. Let $f\in MS(M^n;a,b,1)$ be a Morse-Smale diffeomorphism with a unique saddle $\sigma$. Due to \cite{Robinson99,Smale67}, the separatrices of a saddle of Morse-Smale diffeomorphism have no intersections. Therefore, $f$ is a Morse-Smale diffeomorphism without heteroclinic submanifolds on codimension one separatrices.

It follows from Proposition \ref{prop:rough-decomposition} that $c_1-(a+b)+2\geq 0$ where $c_1\geq 0$ the number of codimension one saddles. Hence, $a+b\leq c_1+2\leq 3$.
Taking in mind that $a\geq 1$ and $b\geq 1$ for any Morse-Smale diffeomorphism, one gets either $a=b=1$ or $a+b=3$ provided $c_1=1$. In the first case $a=b=1$, the result follows from
Proposition \ref{prop:from-Medvedev-Zh-Top-Appl-2013}. In the second case $a+b=3$, one gets $c_1=1=c$ i.e. a unique saddle $\sigma$ is a codimension one saddle, and $c-(a+b)+2 = 0$. It follows that the decomposition (\ref{eq:decomposition-1}) does not hold. Due to Proposition \ref{prop:rough-decomposition}, $M^n$ is a sphere $\mathbb{S}^n$. Clearly that the Morse index of the saddle $\sigma$ equals $(n-1)$ provided $a=1$, $b=2$ or the Morse index of $\sigma$ equals 1 provided $a=2$, $b=1$.
$\Box$

\medskip
\textsl{Proof of Theorem \ref{thm:two-saddles-at-least-one-codimension-one}}. First, we consider case 1) when the saddles $\sigma_1$, $\sigma_2$ of some Morse-Smale diffeomorphism $f\in MS(M^n;a,b,2)$ are both of codimension one. Denote by $Sep_1(\sigma_i)$, $Sep_2(\sigma_i)$ the one-dimensional separatrices of $\sigma_i$, $i=1,2$.
Suppose that
 $$ \left(Sep_1(\sigma_1)\cup Sep_2(\sigma_1)\right)\cap\left(W^s(\sigma_2)\cup W^u(\sigma_2)\right)=\emptyset ,\quad
    \left(Sep_1(\sigma_2)\cup Sep_2(\sigma_2)\right)\cap\left(W^s(\sigma_1)\cup W^u(\sigma_1)\right)=\emptyset $$
As a consequence, there are three possibilities

a) $clos~\left(Sep_1(\sigma_1)\cup Sep_2(\sigma_1)\right)$ and $clos~\left(Sep_1(\sigma_2)\cup Sep_2(\sigma_2)\right)$ are circles;

b) $clos~\left(Sep_1(\sigma_1)\cup Sep_2(\sigma_1)\right)$ and $clos~\left(Sep_1(\sigma_2)\cup Sep_2(\sigma_2)\right)$ are segments;

c) $clos~\left(Sep_1(\sigma_1)\cup Sep_2(\sigma_1)\right)$ is a circle and $clos~\left(Sep_1(\sigma_2)\cup Sep_2(\sigma_2)\right)$ is a segment.

\medskip
In the case (a), we denote the circles $clos\left(Sep_1(\sigma_1)\cup Sep_2(\sigma_1)\right)$ and $clos\left(Sep_1(\sigma_2)\cup Sep_2(\sigma_2)\right)$ by $S_1$ and $S_2$ respectively. Due to
Proposition \ref{prop:closure-is-sphere-for-dim-equals-one}, $S_i$ has a tubular neighborhood $T_i$ homeomorphic to $\mathbb{D}^{n-1}\otimes S^1$, $i=1,2$. Let us prove that one of $T_i$ is forward-invariant while another is back-invariant. Suppose the contrary. Assume for definiteness that the both $T_1$ and $T_2$ are back-invariant i.e. $f^{-1}(T_i)\subset T_i$, $i=1,2$. Note that $f$ has no stable separatrices different from $Sep_1(\sigma_i)$, $Sep_2(\sigma_i)$, $i=1,2$. It follows that there are no sources in $M^n\setminus\left(T_1\cup T_2\right)$ because of the union of all sources and unstable separatrices forms a connected set \cite{GrinesMedvedevPochinkaZh2010}. Since $M^n\setminus\left(S_1\cup S_2\right)$ is connected, the set $M^n\setminus\left(T_1\cup T_2\right)$ contains only one sink. Hence, $f$ has three nodes. On the other side, the unstable manifolds $W^u(\sigma_1)$, $W^u(\sigma_2)$ have no intersections. It follows from Proposition \ref{prop:rough-decomposition} that the number of nodes must be even. This contradiction proves that one of $T_i$, say $T_1$, is forward-invariant while $T_2$ is back-invariant.

Let us prove that there are no fixed points in $M^n\setminus\left(T_1\cup T_2\right)$. Suppose the contrary. Without loss of generality, one can assume that there is a sink
$\omega\in M^n\setminus\left(T_1\cup T_2\right)$. Note that there are sink $\omega\in S_2\subset T_2$ and source $\alpha\in S_1\subset T_1$. Since $W^u(\sigma_2)\subset T_2$, the sink $\omega_0$ belongs to the limit set $Lim\left(W^u(\sigma_1)\right)$ of $W^u(\sigma_1)$. Suppose for a moment that $W^u(\sigma_1)\cap W^s(\sigma_2)\neq\emptyset$. Then
$\omega\subset Lim\left(W^u(\sigma_1)\right)$ that contradicts to a connectedness of the limit set $Lim\left(W^u(\sigma_1)\right)$ of $W^u(\sigma_1)$. Thus,
$W^u(\sigma_1)\cap W^s(\sigma_2)=\emptyset$. It follows that the union $W^u(\sigma_1)\cup\{\omega_0\}=S^{n-1}_0$ is a locally flat embedded $(n-1)$-sphere. This sphere does not divide $M^n$, since
$clos\left(Sep_1(\sigma_1)\cup Sep_2(\sigma_1)\right)$ is a circle. Therefore, the cutting
of $M^n$ along $S^{n-1}_0$ gives a connecter manifold, say $\hat{M}^n$, with two boundary components $M_1$, $M_2$ each homeomorphic to $S^{n-1}_0$. One can glue two $n$-balls $B^n_1$, $B^n_2$ along their boundaries to $M_1$, $M_2$ respectively to get a closed manifold $\tilde{M}^n$. Since $S^{n-1}_0$ is an attracting set, one can extend $f$ on $\tilde{M}^n$ to get a Morse-Smale diffeomorphism $\tilde{f}: \tilde{M}^n\to \tilde{M}^n$ with sinks
$\omega_i\in B^n_i$, $i=1,2$. We see that $\tilde{f}$ has a unique saddle $\sigma_2$, and the source $\alpha$, and at least three sinks $\omega$, $\omega_i$, $i=1,2$. This contradicts to Theorem \ref{thm:one-saddle-structure}.

Thus, $NW(f)=\{\alpha\}\cup\{\omega\}\cup\{\sigma_1\}\cup\{\sigma_2\}\in T_1\cup T_2$. Moreover, $S_1\subset T_1$ is an repelling set while $S_2\subset T_2$ is an attracting set. Since $\partial T_1$ is compact, there is $k\in\mathbb{N}$ such that $f^k(\partial T_1)\subset T_2$. It follows that $f\in MS(M^n;1,1,2)$, and $M^n$ homeomorphic to the union of two copy $\mathbb{D}^{n-1}\otimes S^1$, and the saddles have different Morse indexes.

In the case (b), we denote the segments $clos \left(Sep_1(\sigma_1)\cup Sep_2(\sigma_1)\right)$, $clos \left(Sep_1(\sigma_2)\cup Sep_2(\sigma_2)\right)$ by $I_1$ and $I_2$ respectively. First, we consider the possibility when $I_1$ is an repelling set while $I_2$ is an attracting set. Therefore, there are neighborhoods $U_i$ of $I_i$, $i=1,2$, such that $U_1\subset f(U_1)$ and $f(U_2)\subset U_2$. It follows that one can modify $f$ in $U_1\cup U_2$ to get a Morse-Smale diffeomorphism $\tilde{f}: M^n\to M^n$ such that $\tilde{f}$ has a unique source $\alpha_0\in U_1$, and a unique sink $\omega_0\in U_2$, and $\tilde{f}$ coincides with $f$ out of $U_1\cup U_2$. Since $\tilde{f}$ has no saddles at all, $NW(\tilde{f})=\{\alpha_0\}\cup\{\omega_0\}$. Hence, $\tilde{f}\in MS(\mathbb{S}^n,1,1,0)$. As a consequence, $f\in MS(\mathbb{S}^n,2,2,2)$, and the saddles have different Morse indexes. Now, let us assume that the both $I_1$ and $I_2$ are repelling sets. Taking in mind Theorem \ref{thm:one-saddle-structure}, it follows that $\tilde{f}\in MS(\mathbb{S}^n,1,2,1)\cup MS(\mathbb{S}^n,2,1,1)$. Hence,
$f\in  MS(\mathbb{S}^n,1,3,2)\cup MS(\mathbb{S}^n,3,1,2)$, and the saddles have the same Morse index.

Let us prove that the case (c) does not hold. Suppose the contrary. Applying Proposition \ref{prop:closure-is-sphere-for-dim-equals-one}, one gets a Morse-Smale diffeomorphism
$\tilde{f}\in MS(M^n;1,1,1)$. This contradicts to Theorem \ref{thm:one-saddle-structure}, since the codimension one separatrix (as well as the one-dimensional) of a unique saddle must be $\frac{n}{2}$-dimensional. This is impossible because of $n\geq 4$.

Now, suppose that a codimension one separatrix of one saddle, say $\sigma_1$, intersects a one-dimensional separatrix of $\sigma_2$. It follows that the saddles have the same Morse-index. Without loss of generality, we can assume that the codimension one separatrices of the both saddles $\sigma_1$, $\sigma_2$ are unstable. Since periodic orbits of Morse-Smale diffeomorphism do not form cycles,
$W^u(\sigma_2)\cap W^s(\sigma_1)=\emptyset$. It follows that $clos\left(W^u(\sigma_2)\right)=S^{n-1}$ is a locally flat embedding $(n-1)$-sphere containing a unique sink $\omega$. Cutting $M^n$ along $S^{n-1}$, one gets a manifold $\hat{M}^n$ with two boundary components $M_1$, $M_2$ each homeomorphic to $\mathbb{S}^{n-1}$. One can glue two $n$-balls $B^n_1$, $B^n_2$ along their boundaries to $M_1$, $M_2$ respectively to get a closed manifold $\tilde{M}^n$. Since $S^{n-1}$ is an attracting set, one can extend $f$ on $\tilde{M}^n$ to get a Morse-Smale diffeomorphism $\tilde{f}: \tilde{M}^n\to \tilde{M}^n$ with sinks $\omega_i\in B^n_i$, $i=1,2$ instead of $\sigma_2$ and $\omega$. Note that $\tilde{f}$ has a unique saddle fixed point $\sigma_1$. There are two possibilities: (1) $\tilde{M}^n$ is a connected manifold; (2) $\tilde{M}^n$ consists of two connected manifolds $\tilde{M}^n_1$, $\tilde{M}^n_2$. In the case (1), the non-wandering set $NW(\tilde{f})$ contains a unique saddle $\sigma_1$ with the Morse index $(n-1)$ and at least two sinks $\omega_i$, $i=1,2$. This contradicts to Theorem \ref{thm:one-saddle-structure}. In the case (2), there is component, say $\tilde{M}^n_2$, with no saddle periodic points. Hence, $\tilde{M}^n_2=\mathbb{S}^n$ and the non-wandering set $NW(\tilde{f})\cap \tilde{M}^n_2$ consists of a sink and source. It follows that
$M^n = \tilde{M}^n_1\sharp\mathbb{S}^n$ homeomorphic to $\tilde{M}^n_1$ and $NW(\tilde{f})$ contains a unique sink. Due to Theorem \ref{thm:one-saddle-structure}, $\tilde{f}\in MS(\mathbb{S}^n;1,2,1)$. Hence,
$f\in MS(\mathbb{S}^n,2,2,2)$.

Now, we consider case 2) when one saddle is codimension one while another saddle is not codimension one and two. Assume for definiteness that $\sigma$ is a codimension one saddle while $\sigma_0$ is not codimension one and two.
Suppose also that a codimension one separatrix of $\sigma$ is unstable. Since invariant manifolds of saddles are intersected transversally, $W^s(\sigma)\cap W^u(\sigma_0)=\emptyset$. It follows from
Proposition \ref{prop:closure-is-sphere-for-dim-equals-one} that there are two following possibilities:

2a) $clos(W^s(\sigma))=I$ is a segment;\,\,
2b) $clos(W^s(\sigma))=S^1$ is a circle.

In subcase 2a), the segment $I$ is the union of $W^s(\sigma)$ and two sources, say $\alpha_1$, $\alpha_2$, that are  the endpoints of $I$. By Proposition \ref{prop:closure-is-sphere-for-dim-equals-one}, $I$ has a closed neighborhood $B$ homeomorphic to an $n$-ball such that $int f(B)\subset B$. Therefore, one can modify $f$ inside of $B$ to get a Morse-Smale diffeomorphism $\tilde{f}: M^n\to M^n$ such that $\tilde{f}$ has inside of $B$ a unique source $\alpha_0$. It follows that $\tilde{f}$ has a unique saddle $\sigma_0$ that is not a saddle of codimension one. Due to Proposition \ref{prop:from-Medvedev-Zh-Top-Appl-2013}, $n\in\{4,8,16\}$, and $M^n$ is a disjoint union of an open ball $\mathbb{B}^n$ and sphere $\mathbb{S}^{\frac{n}{2}}$. In addition, the homotopy groups $\pi_1(M^n)=\cdots=\pi_{\frac{n}{2}-1}(M^n)=0$, and hence, $M^n$ is simply connected and orientable. Moreover,
$\tilde{f}\in MS(M^n;1,1,1)$. Hence, $f\in MS(M^n;1,2,2)$. Clearly, if we assume that a codimension one separatrix of $\sigma$ is stable, one gets $f\in MS(M^n;2,1,2)$.

Let us show that the subcase 2b) does not hold. Suppose the contrary. The circle $S^1$ is the union of $W^s(\sigma)$ and a source, say $\alpha_0$. Since $\alpha_0$ is a hyperbolic source, there is a neighborhood $B_0$ of $\alpha_0$ such that $B$ is an $n$-ball and $clos f^{-1}(B)\subset B$. Without loss of generality, one can assume that the boundary $\partial B$ is a smoothly embedded $(n-1)$-sphere. Since $n\geq 4$, $S^1$ is a locally flat embedded circle. Therefore, one can slightly modify $\partial B$ to get the transversal intersection $\partial B\pitchfork W^s(\sigma)$ consisting of two points, say $b_1$ and $b_2$. The intersection $W^s(\sigma)\cap B_0$ is the path denoted by $C$ connecting the points $b_1$ and $b_2$. According to the Grobman-Hartman theorem, $f$ conjugate to a linear hyperbolic mapping in a sufficiently small neighborhood of $\alpha_0$. Therefore, there is a tube $T$ homeomorphic $C\times d^{n-1}$ such that $T\subset clos B_0$, and $f^{-1}(T)\subset int T$ where $d^{n-1}$ is a closed $(n-1)$-disk. The tube $T$ can be continued to a closed neighborhood $N$ of $S^1$ such that $N$ homeomorphic $N\times d^{n-1}$ such that $f^{-1}(N)\subset int N$. It follows that $S^1$ is a repelling invariant set with the basin $B(S^1)=\cup_{m\geq 0}f^m(N)$. By construction, $S^1$ does not homotopy to zero in $B(S^1)$.

Now, we'll prove that $S^1$ homotopy to zero in $B(S^1)$. First, let us show that $f$ has a unique sink $\omega_0$ such that
\begin{equation}\label{eq:1}
    M^n\setminus S^1\subset\left(W^s(\omega_0)\cup W^s(\sigma_0)\right).
\end{equation}
Suppose for a moment that there exists one more sink $\omega_1\neq\omega_0$. The union $G=W^u(\sigma_1)\cup W^u(\sigma_0)\cup\{\omega_1\}\cup\{\omega_0\}$ is a connected global attractor of $f$ \cite{GrinesMedvedevPochinkaZh2010}.
Because of the unstable manifold $W^u(\sigma_0)$ has no heteroclinic intersection, $W^u(\sigma_0)\cup\{\omega_0\}$ is a topologically embedded sphere such that $W^u(\sigma_0)\setminus\{\omega_0\}\subset W^s(\omega_0)$. Thus,
\begin{equation}\label{eq:2}
    W^u(\sigma_0)\cap W^s(\omega_0)=\emptyset.
\end{equation}
The connectedness of $G_0$ implies that $W^u(\sigma)\cap W^s(\omega_0)\neq\emptyset$ and $W^u(\sigma)\cap W^s(\omega_1)\neq\emptyset$. This means that there is a saddle $\sigma_*$ such that
$\{\omega_1\}\cup\{\omega_0\}\subset clos W^u(\sigma_*)$, and $W^u(\sigma)\cap W^s(\omega_*)\neq\emptyset$. Hence, $\sigma_*=\sigma_0$ provided there are only two saddles $\sigma_0$ and $\sigma_1$. The inclusion
$\{\omega_1\}\cup\{\omega_0\}\subset clos W^u(\sigma_*)$ means that $W^u(\sigma_0)\cap W^s(\omega_0)\neq\emptyset$ that contradicts to (\ref{eq:2}). We see that $\omega_0$ is a unique saddle of $f$. Since $M^n$ is the disjoint union of non-wandering points, one get (\ref{eq:1}).

It follows from (\ref{eq:1}) that the basin $B(S^1)$ consists of $S^1$, and $W^s(\omega_0)\setminus W^u(\sigma_0)$, and the intersection $W^u(\sigma)\cap W^s(\sigma_0)$ if nonempty. The dimension of the stable manifold $W^s(\sigma_0)$ is not more than $n-2$. Therefore, there is a circle $\tilde{S}^1$ belonging to the open domain $W^s(\omega_0)\setminus W^u(\sigma_0)$ such that $\tilde{S}^1$ homotopy $S^1$. One can choose $\tilde{S}^1$ being arbitrary close to $S^1$ so that $\tilde{S}^1\subset B(S^1)$. Since $\{\omega_0\}\cup W^u(\sigma_0)$ is a topologically embedded sphere whose dimension no more than $n-2$, it is possible to deform $\tilde{S}^1$ to $W^s(\omega_0)$ so that $\tilde{S}^1\subset B(S^1)\cap W^s(\omega_0)$. We know that $W^s(\omega_0)$ is homeomorphic to $\mathbb{R}^n$. Hence, there is an immersion $\psi: D^2\to W^s(\omega_0)$ such that $\psi|_{\partial D^2}=\tilde{S}^1$ where $D^2$ is a closed 2-disk. Slightly moving $\psi(D^2)$, one can assume that $\psi(D^2)$ is transversal to $W^s(\sigma_0)$ \cite{Haefliger1961a}. Since $\dim W^s(\sigma_0)\geq n-3$, the transversality of the intersection $\psi(D^2)\cap W^s(\sigma_0)$ means that $\psi(D^2)\cap W^s(\sigma_0)=\emptyset$. It follows that $\tilde{S}^1$ homotopy to zero in $B(S^1)$. Hence, $S^1$ homotopy to zero in $B(S^1)$ as well. This contradiction concludes the proof.
$\Box$

\medskip
\textsl{Proof of Theorem \ref{thm-all-saddles-codimension-one-same-morse-index}}. From Theorems \ref{thm:one-saddle-structure}, \ref{thm:two-saddles-at-least-one-codimension-one} it follows that the desired result holds for $c=1,2$. Assume that the result holds for $c=1,2,\ldots,k$. Take $f\in MS(M^n;a,b,k+1)$ with codimension one saddles $\sigma_1$, $\sigma_2$, $\ldots$, $\sigma_{k+1}$. For definiteness, assume that the Morse index of each
$\sigma_i$, $1\leq i\leq k+1$, equals $n-1$. Hence, there is a saddle, say $\sigma_{k+1}$, such that the unstable manifold $W^u(\sigma_{k+1})$ has no heteroclinic intersections. According
Proposition \ref{prop:closure-is-sphere-for-dim-more-one}, the topological closure $clos W^u(\sigma_{k+1})=S^{n-1}$ is a locally flat embedded $(n-1)$-sphere. Cutting $M^n$ along $S^{n-1}$, one gets a manifold $\hat{M}^n$ with two boundary components $M_1$, $M_2$ each homeomorphic to $\mathbb{S}^{n-1}$. One can glue two $n$-balls $B^n_1$, $B^n_2$ along their boundaries to $M_1$, $M_2$ respectively to get a closed manifold $\tilde{M}^n$. Since $S^{n-1}$ is an attracting set, one can extend $f$ on $\tilde{M}^n$ to get a Morse-Smale diffeomorphism $\tilde{f}: \tilde{M}^n\to \tilde{M}^n$ with sinks $\omega_i\in B^n_i$, $i=1,2$ instead of $\sigma_{k+1}$ and some sink. There are two possibilities: (1) $\tilde{M}^n$ is a connected manifold; (2) $\tilde{M}^n$ consists of two connected manifolds $\tilde{M}^n_1$, $\tilde{M}^n_2$.

In the case (1), the non-wandering set $NW(\tilde{f})$ contains $k$ codimension one saddles $\sigma_1$, $\ldots$, $\sigma_k$ with the Morse index $n-1$ and at least two sinks $\omega_i\in B^n_i$, $i=1,2$. However, by the inductive assumption, $\tilde{f}$ has a unique sink. This contradiction shows that the case (1) is impossible.

In the case (2), each $\tilde{f}|_{\tilde{M}^n_1}: \tilde{M}^n_1\to\tilde{M}^n_1$, $\tilde{f}|_{\tilde{M}^n_2}: \tilde{M}^n_2\to\tilde{M}^n_2$ has $\leq k$ saddles. By the inductive assumption, $\tilde{M}^n_1=\tilde{M}^n_1=\mathbb{S}^n$. Hence, $M^n=\mathbb{S}^n\sharp\mathbb{S}^n=\mathbb{S}^n$. Moreover, $\tilde{f}|_{\tilde{M}^n_1}\in MS(\mathbb{S}^n;1,b_1,k_1)$ and
$\tilde{f}|_{\tilde{M}^n_2}\in MS(\mathbb{S}^n;1,b_2,k_2)$ where $b=b_1+b_2$, $k_1+k_2=k$, and $b_i=k_i+1$, $i=1,2$. We see that $b=k_1+k_2+2=k+1$. Thus, $f\in MS(\mathbb{S}^n,1,k+1,k)$. This completes the proof.
$\Box$

\medskip
\textsl{Proof of Theorem \ref{thm:there-are-no-cod1-at-all}}. The qualities $\mu=\nu=1$ was proved in \cite{GrinesMedvedevPochinkaZh2010}.
Let $e: S^1\to M^n$ be a map representing an element of $\pi_1(M^n)$. Since $n\geq 4$, one can assume that $e$ is a smooth immersion \cite{Haefliger1961a}. Without loss of generality, one can assume that $e(S^1)$ does not contain periodic points of $f$. Take an unstable manifold $W^u(\sigma_1)$ of some saddle $\sigma_1$. Since $W^u(\sigma_1)$ is an image of $\mathbb{R}^k$ under a smooth immersion, one can slightly move $e(S^1)$ to become transversal to $W^u(\sigma_1)$. By condition, $\dim W^u(\sigma_1)\leq n-2$. The transversality gives that
$e(S^1)\cap W^u(\sigma_1)=\emptyset$. Continuing this procedure, one can obtain $e(S^1)$ with no intersections with unstable manifolds of every saddles. Hence, $e(S^1)$ belongs to the unstable manifold $W^u(\alpha)$ of unique source $\alpha$ of $f$. Since $W^u(\alpha)$ homeomorphic to $\mathbb{R}^n$, $e(S^1)$ homotopy to zero. It follows that $\pi_1(M^n)=0$.
$\Box$

\section{Examples}\label{s:exm-and-conclusions}

Here, we briefly describe some examples to show that the sets of Morse-Smale diffeomorphisms under consideration are not empty.

\begin{figure}[h]
\centerline{\includegraphics[height=4cm]{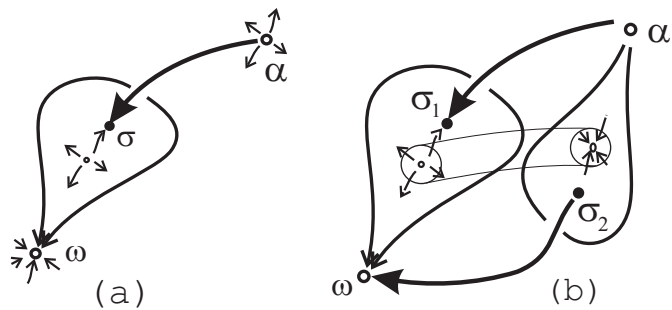}}\caption{(a) $MS(\mathbb{S}^n;1,2,1)\neq\emptyset$, $MS(\mathbb{S}^n;2,1,1)\neq\emptyset$; (b) $MS(\mathbb{S}^{n-1}\otimes S^1;2,2,2)\neq\emptyset$.} \label{1+2-saddles}
\end{figure}

1) $MS(\mathbb{S}^n;1,2,1)\neq\emptyset$, $MS(\mathbb{S}^n;2,1,1)\neq\emptyset$, $MS(M^n;1,1,1)\neq\emptyset$.

On Fig.~\ref{1+2-saddles},(a), one represents the phase portrait of $f\in MS(\mathbb{S}^n;2,1,1)$ with a unique codimen\-sion one saddle $\sigma$, and a sink $\omega$, and two sources. Thus, $MS(\mathbb{S}^n;2,1,1)\neq\emptyset$. Clearly, $f^{-1}\in MS(\mathbb{S}^n;1,2,1)$, and so $MS(\mathbb{S}^n;1,2,1)\neq\emptyset$. Note that it holds for $n=2,3$ as well.
The existence of Morse-Smale diffeomorphisms $f\in MS(M^n;1,1,1)$ was proved in \cite{MedvedevZhuzhoma2013-top-appl}.

2) $MS(\mathbb{S}^n;2,2,2)\neq\emptyset$, $MS(\mathbb{S}^n;1,3,2)\neq\emptyset$, $MS(\mathbb{S}^n;3,1,2)\neq\emptyset$, $MS(M^n;1,2,2)\neq\emptyset$, $MS(M^n;2,1,2)\neq\emptyset$.

\begin{figure}[h]
\centerline{\includegraphics[height=4cm]{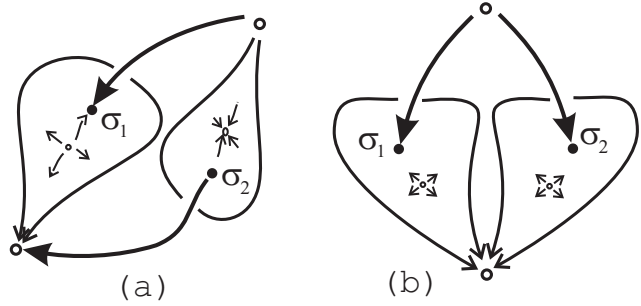}}\caption{(a) $MS(\mathbb{S}^n;2,2,2)\neq\emptyset$; (b) $MS(\mathbb{S}^n;1,3,2)\neq\emptyset$, $MS(\mathbb{S}^n;3,1,2)\neq\emptyset$} \label{S-n}
\end{figure}

On Fig.~\ref{S-n},(a), one illustrates the phase portrait of $f\in MS(\mathbb{S}^n;2,2,2)$ with codimension one saddles $\sigma_1$, $\sigma_2$, and two sinks, and two sources. Hence, $MS(\mathbb{S}^n;2,2,2)\neq\emptyset$. On Fig.~\ref{S-n},(b), one represents $f\in MS(\mathbb{S}^n;1,3,2)$. Clearly, $f^{-1}\in MS(\mathbb{S}^n;3,1,2)$. Thus, $MS(\mathbb{S}^n;1,3,2)\neq\emptyset$ and $MS(\mathbb{S}^n;3,1,2)\neq\emptyset$.

Suppose $f_1\in MS(M^n;1,1,1)$ satisfies Theorem \ref{thm:one-saddle-structure}, item 1), and $f_2\in MS(\mathbb{S}^n;1,2,1)$. Since $M^n=M^n\sharp\mathbb{S}^n$, one can construct $f=f_1\sharp f_2\in MS(M^n;1,2,2)$ satisfying Theorem \ref{thm:two-saddles-at-least-one-codimension-one}, item 2). Thus, $MS(M^n;1,2,2)\neq\emptyset$. Similarly, $MS(M^n;2,1,2)\neq\emptyset$.

3) $MS(\mathbb{D}^{n-1}\otimes S^1\bigcup\mathbb{D}^{n-1}\otimes S^1;1,1,2)\neq\emptyset$.

\begin{figure}[h]
\centerline{\includegraphics[height=4cm]{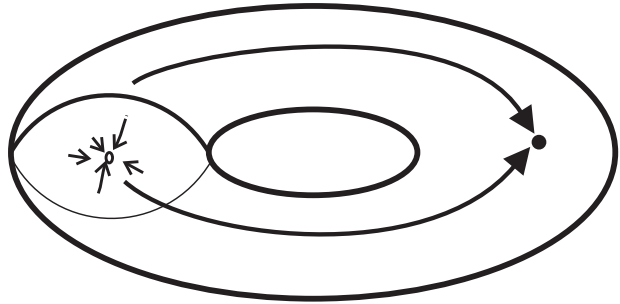}}\caption{$MS(\mathbb{D}^{n-1}\otimes S^1\bigcup\mathbb{D}^{n-1}\otimes S^1;1,1,2)\neq\emptyset$.} \label{solid-tor}
\end{figure}

Take a high-dimensional solid torus $T_1=\mathbb{D}^{n-1}\otimes S^1$ supporting a vector field $V_1$ with the phase portrait represented by Fig.~\ref{solid-tor}. Let $T_2$ be a copy of $T_1$ with the vector field $-V_1$. Gluing carefully $T_1$, $T_2$ along their boundaries one gets $M=\mathbb{D}^{n-1}\otimes S^1\bigcup\mathbb{D}^{n-1}\otimes S^1$ the Morse-Smale vector field denoted by $V$. The time-shift along the trajectories of $V$ is a Morse-Smale diffeomorphism belonging to $MS(\mathbb{D}^{n-1}\otimes S^1\bigcup\mathbb{D}^{n-1}\otimes S^1;1,1,2)$. Another way to construct particular cases
$MS(\mathbb{S}^{n-1}\otimes S^1;1,2,2)\subset MS(\mathbb{D}^{n-1}\otimes S^1\bigcup\mathbb{D}^{n-1}\otimes S^1;1,1,2)$ is illustrated in Fig.~\ref{1+2-saddles},(b).

4) $MS(\mathbb{S}^n;c+1,1,c)\neq\emptyset$, $MS(\mathbb{S}^n;1,c+1,c)\neq\emptyset$.

Starting with $f_1\in MS(\mathbb{S}^n;1,3,2)$, see Fig.~\ref{S-n},(b), and adding $f_2\in MS(\mathbb{S}^n;1,2,1)$, see Fig.~\ref{1+2-saddles},(a), one can get $f\in MS(\mathbb{S}^n;1,4,3)$. Continuing this way, we get that $MS(\mathbb{S}^n;1,c+1,c)\neq\emptyset$ for any $c\geq 3$. Similarly, $MS(\mathbb{S}^n;c+1,1,c)\neq\emptyset$.

%%%%%%%%%%%%%%%%%%%%%%%%%%%%%%%%%%%%%%%%%%%%%%%%%%%%%%%%%%%%%%%%%%%%%%%%%%%

\bigskip
\noindent
National Research University Higher School of Economics, 25/12 Bolshaya Pecherskaya,\\ 603005, Nizhni Novgorod, Russia

\bigskip
\noindent
\textit{E-mails:} medvedev-1942@mail.ru, zhuzhoma@mail.ru % medvedev@uic.nnov.ru

\end{document}